\documentclass{article} 
\usepackage{nips15submit_e,times}
\usepackage{hyperref}
\usepackage{url}

\usepackage{graphicx} 
\usepackage{subfigure}

\usepackage{natbib}

\usepackage{algorithm}
\usepackage{algorithmic}
\usepackage{amsmath, amsfonts, amssymb, amsthm}
\usepackage[misc,geometry]{ifsym}

\def\R{\mathbb R}
\newcommand{\E}{{\mathbb E}}
\def\reff#1{{\rm (\ref{#1})}}  
\def\myproof{{\bf Proof.\ }}
\def\e{\varepsilon}
\def\vp{\varphi}
\def\tpi{\tilde{\pi}}
\def\B{{\mathcal B}}
\def\S{{\mathcal S}}
\def\U{\mathcal{U}}
\newcommand{\la}{\langle}
\newcommand{\ra}{\rangle}

\numberwithin {equation}{section}

\newtheorem{Lm}{Lemma}

\newtheorem{Th}{Theorem}

\nipsfinalcopy

\title{Learning Supervised PageRank with Gradient-Free Optimization Methods}

\author{
Lev Bogolubsky \\
Yandex\\
Leo Tolstoy st. 16\\
Moscow, Russian Federation\\
\texttt{bogolubsky@yandex-team.ru} \\
\And
Pavel Dvurechensky
\\
Weierstrass Institute for Applied Analysis and Stochastics\\
Mohrenstr. 39, 10117 Berlin, Germany; \\
Institute for Information Transmission Problems RAS\\
Bolshoy Karetny per. 19, build.1\\
Moscow 127051 Russia \\
\texttt{pavel.dvurechensky@wias-berlin.de} \\
\And
Alexander Gasnikov \\
Institute for Information Transmission Problems RAS \\
Bolshoy Karetny per. 19, build.1\\
Moscow 127051 Russia \\
\texttt{gasnikov@yandex.ru} \\
\AND
Gleb Gusev \\
Yandex\\
Leo Tolstoy st. 16\\
Moscow, Russian Federation\\
\texttt{gleb57@yandex-team.ru} \\
\And
Yurii Nesterov\\
Center for Operations Research and Econometrics (CORE) \\
34 voie du Roman Pays, 1348, Louvain-la-Neuve, Belgium \\
\texttt{yurii.nesterov@uclouvain.be} \\
\And
Andrey Raigorodskii\\
Yandex\\
Leo Tolstoy st. 16\\
Moscow, Russian Federation\\
\texttt{raigorodsky@yandex-team.ru} \\
\And
Aleksey Tikhonov\\
Yandex\\
Leo Tolstoy st. 16\\
Moscow, Russian Federation\\
\texttt{altsoph@yandex-team.ru}\\
\And
Maksim Zhukovskii \\
Yandex\\
Leo Tolstoy st. 16\\
Moscow, Russian Federation\\
\texttt{zhukmax@yandex-team.ru} \\
}

\nipsfinalcopy 

\begin{document}

\maketitle

\begin{abstract}
In this paper, we consider a problem of learning supervised PageRank models, which can account for some properties not considered by classical approaches such as the classical PageRank algorithm. Due to huge hidden dimension of the optimization problem we use random gradient-free methods to solve it. We prove a convergence theorem and estimate the number of arithmetic operations needed to solve it with a given accuracy.  We find the best settings of the gradient-free optimization method in terms of the number of arithmetic operations needed to achieve given accuracy of the objective. In the paper, we apply our algorithm to the web page ranking problem. We consider a parametric graph model of users' behavior and evaluate web pages' relevance to queries by our algorithm. The experiments show that our optimization method outperforms the untuned gradient-free method in the ranking quality.
\end{abstract}

\section{Introduction}
\label{intro}

The most acknowledged methods of measuring importance of nodes in graphs are based on random walk models. Particularly, PageRank~[18], HITS~[11], and their variants~[8,~9,~19] are originally based on a discrete-time Markov random walk on a link graph. 
According to the PageRank algorithm, the score of a node equals to its probability in the stationary distribution of a Markov process, which models a random walk on the graph. 
Despite undeniable advantages of PageRank and its mentioned modifications, 
these algorithms miss important aspects of the graph that are not described by its structure. 

In contrast, a number of approaches allows to account for different properties of nodes and edges between them by encoding them in restart and transition probabilities (see~[3,~4,~6,~10,~12,~20,~21]). These properties may include, e.g., the statistics about users' interactions with the nodes (in web graphs~[12] or graphs of social networks~[2]), types of edges (such as URL redirecting in web graphs~[20]) or histories of nodes' and edges' changes~[22]. 
Particularly, the transition probabilities in BrowseRank algorithm~[12] are proportional to weights of edges which are equal to numbers of users' transitions.
In the general ranking framework called Supervised PageRank~[21], weights of nodes and edges in a graph are linear combinations of their features with coefficients as the model parameters. 
The authors consider an optimization problem for learning the parameters and solve it by a gradient-based optimization method. 
However, this method is based on computation of derivatives of stationary distribution vectors w.r.t. its parameters which include calculating the derivative for each element of a billion by billion matrix and, therefore, seems to be computationally very expensive. The same problem appears when using coordinate descent methods like~[15] does. Another obstacle to the use of gradient or coordinate descent methods is that we can't calculate derivatives precisely, since we can't evaluate the exact stationary distribution.

In our paper, we consider the optimization problem from~[21] and propose a two-level method to solve it. On the lower level, we use the linearly convergent method from~[17] to calculate an approximation to the stationary distribution of the Markov process. We show in Section~\ref{learn} that this method has the best among others~[5] complexity bound for the two-level method as a whole. However, it is not enough to calculate the stationary distribution itself, since we need also to optimize the parameters of the random walk with respect to an objective function, which is based on the stationary distribution. To overcome the above obstacles, we use a gradient-free optimization method on the upper level of our algorithm. The standard gradient-free optimization methods~[7,~16] require exact values of the objective function. 
Our first contribution described in Section~\ref{method} consists in adapting the framework of~[16] to the case when the value of the function is calculated with some known accuracy. We prove a convergence theorem (Section~\ref{method}) for this method. Our second contribution consists in investigating the trade-off between the accuracy of the lower level algorithm, which is controlled by the number of iterations, and the computational complexity of the two-level algorithm as a whole (Section~\ref{learn}). For given accuracy, we estimate the number of arithmetic operations needed by our algorithm to find the values of parameters such that the difference between the respective value of the objective and its local minimum does not exceed this accuracy. 
In the experiments, we apply our algorithm to the problem of web pages' ranking. 
We show in Section~\ref{results} that our two-level method outperforms an untuned gradient-free method in the ranking quality. 

The remainder of the paper is organized as follows. In Section~\ref{model}, we describe the random walk model. In Section~\ref{optimal}, we define the learning problem and discuss its properties and possible methods for its solution. In Section~\ref{method} we describe the framework of random gradient-free optimization methods and generalize it to the case when the function values are inaccurate.
In Section~\ref{learn} we propose two-level algorithm for the stated learning problem.
The experimental results are reported in Section~\ref{Experimental results}. In Section~\ref{Conclusion}, we summarize the outcomes of our study,
discuss its potential applications and directions of future work.

\section{Model description}
\label{model}

Let $\Gamma=(V,E)$ be a directed graph. 
Denote by $p$ the number of vertices in $V$. 
Let
$$
\mathcal{F}_1=\{F(\varphi_1,\cdot): V \rightarrow\mathbb{R}\},\,\,\,
\mathcal{F}_2=\{G(\varphi_2,\cdot): E \rightarrow\mathbb{R}\}
$$
be two classes of functions parameterized by $\vp_1\in\mathbb{R}^{m_1}, \vp_2\in\mathbb{R}^{m_2}$ respectively, where $m_1$ is the number of nodes' features, $m_2$ is the number of edges' features. We denote $m=m_1+m_2$ 
and $\varphi = (\vp_1, \vp_2)^T$. Let us describe the random walk on the graph $\Gamma$, which was considered in~[21]. The {\it seed set} $V^1\subset V$ is defined as follows: $i\in V^1$ if and only if $F(\varphi_1,i)\neq 0$ for some $\varphi_1\in\mathbb{R}^{m_1}$. A surfer starts a random walk from a random page $i\in V^1$, the initial probability of being at vertex $i$ is called the {\it restart probability} and equals
\begin{equation}
 [\pi^0(\vp)]_i=\frac{F(\vp_1, i)}{\sum_{\tilde i\in  V^1}F(\vp_1, \tilde i)}
\label{restart}
\end{equation}
(equals $0$ for $i\in V\setminus V^1$). 
At each step, the surfer (with a current position $\tilde i\in V$) either chooses any vertex from $V^1$ in accordance with the distribution $\pi^0(\vp)$ (makes a {\it restart}) with probability $\alpha\in(0,1)$, which is called the {\it damping factor}, or chooses to traverse an outcoming edge (makes a {\it transition}) with probability $1-\alpha$. 
The probability
\begin{equation}
 [P(\vp)]_{\tilde i,i}=\frac{G(\vp_2, \tilde{i} \to i)}{\sum_{j: \tilde{i} \to  j}G(\vp_2,\tilde{i} \to j)}
\label{transition}
\end{equation}
of traversing an edge $\tilde i\rightarrow i\in E$ is called the transition probability. Finally, by Equation~\ref{restart} and Equation~\ref{transition} the total probability of choosing vertex $i\in V^1$ conditioned by the surfer being at vertex $\tilde i$ equals $\alpha[\pi^0(\vp)]_i+(1-\alpha)[P(\vp)]_{\tilde i, i}$
 (originally~[18], $\alpha=0.15$).
 If $i\in V\setminus V^1$, then this probability equals $(1-\alpha)[P(\vp)]_{\tilde i,i}$.
Denote by $\pi \in \R^{p}$ the stationary distribution of the described Markov process.
It can be found as a solution of the system of equations
\begin{equation}
[\pi]_i = \alpha [\pi^0(\vp)]_i + (1-\alpha)\sum_{\tilde{i}: \tilde{i} \to i \in E}[P(\vp)]_{\tilde i,i}
[\pi]_{\tilde i}.
\label{pi_general}
\end{equation}
In this paper, we learn the ranking algorithm, which orders the vertices $i$ by their probabilities $[\pi]_{i}$ in the stationary distribution $\pi$.


\section{Learning problem statement}
\label{optimal}

Let $Q$ be a set of search queries and weights of nodes and edges $F_q:=F$ and $G_q:=G$ depend on $q\in Q$. Let $V_q$ be a set of vertices which are relevant to $q$. In other words, for any $i\in V_q$ either $F_q(\vp_1,i)\neq 0$ for some $\varphi_1\in\mathbb{R}^{m_1}$ or there exists a path $i_0\rightarrow i_1,\ldots, i_k\rightarrow i_{k+1}=i$ in $\Gamma$ such that $F_q(\vp_1,i_0)\neq 0$, $G_q(\vp_2,i_j\rightarrow i_{j+1})\neq 0$ for some $\vp\in\mathbb{R}^{m}$ and all $j\in\{0,\ldots,k\}$. Denote $E_q$ a set of all edges $\tilde i\rightarrow i$ from $E$ such that $\tilde i,i\in V_q$ and $G_q(\varphi_2,\tilde i\rightarrow i)\neq 0$ for some $\varphi_2\in\mathbb{R}^{m_2}$. For any $q\in Q$, denote $\Gamma_q=(V_q,E_q)$. For fixed $q\in Q$, the graph $\Gamma_q$ and functions $F_q,G_q$, we consider the notations from the previous section and add the index $q$: $V_q^1:=V^1$, $\pi_q^0:=\pi^0$, $P_q:=P$, $p_q:=p$, $\pi_q:=\pi$. The parameters $\alpha$ and $\vp$ of the model do not depend on $q$.

Our goal is to find the  parameters vector $\vp$ which minimizes the discrepancy of the nodes ranking scores $[\pi_q]_i$, $i \in V_q$, calculated as the stationary distribution in the above Markov process from the nodes ranking scores defined by assessors. For each $q\in Q$, there is a set of nodes in $V_q$ manually judged and grouped by relevance labels $1,\ldots,k$. We denote $V^j_q$ the set of documents annotated with label $k+1-j$ (i.e., $V_q^1$ is the set of all nodes with the highest relevance score). For any two nodes $i_1\in V_q^{j_1},i_2\in V_q^{j_2}$, let $h(j_1,j_2,[\pi_q]_{i_2}-[\pi_q]_{i_1})$ be the value of the loss function. If it is non-zero, then the position of the node $i_1$ according to our ranking algorithm is higher than the position of the node $i_2$ but $j_1>j_2$. We consider square loss with margins $b_{j_1j_2}>0$, where $1\leq j_2<j_1\leq k$:
$h(j_1,j_2,x)=(\min\{x+b_{j_1j_2},0\})^2$ as it was done in previous studies~[12,~21,~22]. We minimize
\begin{equation}
 f(\vp)= \frac{1}{|Q|} \sum_{q=1}^{|Q|}\sum\limits_{1\leq j_2<j_1\leq k}\sum\limits_{i_1\in V_q^{j_1},i_2\in V_q^{j_2}}h(j_1,j_2,[\pi_q]_{i_2}-[\pi_q]_{i_1})
 \label{eq:f_phi_def_1}
\end{equation}
in order to learn our model using the data given by assessors. 

As it was said above, finding nodes ranking scores for the fixed query $q$ leads to the problem of finding the stationary distribution $\pi_q$ of the Markov process as a solution of Equation~\ref{pi_general}
or equivalently
\begin{equation}
\pi_q = \alpha \pi^0_q(\vp)  + (1-\alpha) P_q^T(\varphi)  \pi_q.
\label{eq:pi_Phi_P2}
\end{equation}
The solution $\pi_q(\vp) $ of \reff{eq:pi_Phi_P2} can be found as
$\pi_q (\vp) = \alpha \left[I - (1-\alpha) P_q^T(\varphi) \right]^{-1} \pi^0_q(\vp)$,
where $I$ is the identity matrix.

It is easy to show~[17] that the vector
\begin{equation}
\tpi^N_q (\vp) = \frac{\alpha}{1-(1-\alpha)^{N+1}} \sum_{i=0}^{N} {(1-\alpha)^i \left[ P_q^T(\varphi) \right]^i \pi^0_q(\vp)}
\label{eq:tpi_def}
\end{equation}
satisfies $\| \tpi^N_q (\vp) - \pi_q (\vp) \|_1 \leq 2 (1-\alpha)^{N+1}$.
As it also shown there to obtain vector $\tpi^N_q (\vp)$ satisfying  
\begin{equation}
\| \tpi^N_q (\vp) - \pi_q (\vp)\|_1 \leq \Delta
\label{Delta}
\end{equation}
one needs $\frac{1}{\alpha} \ln \frac{2}{\Delta}$ iterations of simple iteration method. Each iteration of such method requires one multiplication of the matrix $P_q^T(\varphi)$ by the vector of dimension $p_q$. This requires $s_qp_q$ arithmetic operations. Here $s_q$ is the maximum number of non-zero elements over columns of the matrix $P_q(\vp)$ (the {\it sparsity parameter}).  So the total number of arithmetic operations for obtaining approximation satisfying \reff{Delta} is $\frac{s_qp_q}{\alpha} \ln \frac{2}{\Delta}$ arithmetic operations.  Note that $ s_q \ll p_q$ and that this algorithm for finding the vector $\tpi^N_q (\vp)$ can be fully paralleled.

Let us now turn to the problem of the minimization of the function $f(\varphi)$ \reff{eq:f_phi_def_1}. We can rewrite this function as
\begin{equation}
f(\varphi)= \frac{1}{|Q|} \sum_{q=1}^{|Q|} \|(A_q \pi_q (\vp) +b_q)_{+} \|^2_2,
\label{eq:f_phi_def_2}
\end{equation}
where vector $x_+$ has components $[x_+]_i=\max\{x_i,0\}$, the matrix $A_q \in \R^{r_q \times p_q}$ represents assessor's view of the relevance of pages to the query $q$, vector $b_q$ is the vector composed from thresholds $b_{j_1,j_2}$ in \reff{eq:f_phi_def_1} with fixed $q$, $r_q$ is the number of summands in \reff{eq:f_phi_def_1} with fixed $q$. 



Due to huge hidden dimension $p_q$, the calculation of the  of $f(\varphi)$ includes calculating the derivative for each element of the $p_q\times p_q$ matrix $P_q(\vp)$ which is too expensive. So we are going to use gradient-free methods for minimization of the function $f(\varphi)$. Such methods were introduced rather long ago, see, e.g.,~[13]. Note that we have to work in the framework of non-exact zero-order oracle. Note that each row of the matrix $A_q$ contains one $1$ and one $-1$, and all other elements of the row are equal to $0$ and hence $\|A_q\|_2 \leq \sqrt{2 r_q}$. This leads to the following Lemma which says how the error of the approximation of $\pi_q(\vp)$ affects the error in the value of the function $f(\varphi)$.

\begin{Lm}
Assume that the vector $\tpi^N_q (\vp)$ satisfies Equation~\ref{Delta}. Denote $r= \max_{q} r_q$, $b = \max_q \|b_q\|_2$. Then
\begin{equation}
f^{\delta}(\varphi) = \frac{1}{|Q|} \sum_{q=1}^{|Q|} \|(A_q \tpi^N_q (\vp) +b_q)_{+} \|^2_2
\label{eq:fdvpDef}
\end{equation}
satisfies $|f^{\delta}(\varphi)  - f(\varphi)| \leq \delta = \Delta \sqrt{2 r} (2 \sqrt{2 r}  + 2 b)$.
\label{Lm:Delta_to_delta}
\end{Lm}

\section{Random gradient-free optimization methods}
\label{method}

Let us describe the well-known framework of random gradient-free methods~[1,~7,~16]. Our main contribution, described in this section, consists in developing this framework for the situation of presence of error of unknown nature in the objective function value. Apart from~[16] we consider randomization on a Euclidean ball which seems to give better large deviations bounds and doesn't need the assumption that the function can be calculated at any point of the space $\R^m$.

In this section, we consider a general function $f(\cdot)$ and denote its argument by $x$ or $y$ to avoid confusion with other sections.
Assume that the function $f (\cdot): \R^m \to \R$ is convex and has Lipschitz continuous gradient with constant $L$ (we write $f \in C^{1,1}_{L}$):
\begin{equation}
     |f(x)-f(y) - \la \nabla f(y) ,x-y \ra| \leq \frac{L}{2}\|x-y\|_2^2, \quad x,y \in \R^m.
     \notag
     \label{eq:fLipSm}
\end{equation}
Also we assume that the oracle returns the value $f^{\delta}(x)=f(x) + \tilde{\delta}(x)$, where $\tilde{\delta}(x)$ is the oracle error satisfying $|\tilde{\delta}(x)| \leq \delta$.
Consider smoothed counterpart of the function $f(x)$:
\begin{equation}
f_{\mu}(x) = \E f(x+\mu\xi) = \frac{1}{V_{\B}} \int_{\B} f(x+\mu t) dt,
\notag
\end{equation}
where $\xi$ is uniformly distributed over unit ball $\B=\{ t \in \R^m  : \|t\|_2 \leq 1\}$ random vector, $V_{\B}$ is the volume of the unit ball $\B$, $\mu \geq 0$ is a smoothing parameter.
It is easy to show that
\begin{itemize}
     \item If $f$ is convex, then $f_{\mu}$ is also convex
     \item If $f \in C^{1,1}_{L}$, then $f_{\mu} \in C^{1,1}_{L}$.
                 \item If $f \in C^{1,1}_{L}$, then $ f(x) \leq f_{\mu}(x) \leq f(x) + \frac{L \mu^2}{2}$ for all $x \in \R^m$.
\end{itemize}

The random gradient-free oracle is defined as follows
\begin{equation}
g_{\mu}(x) = \frac{m}{\mu}(f(x+\mu s) -f(x)) s,
\notag
\end{equation}
where $s$ is uniformly distributed vector over the unit sphere $\S =\{ t \in \R^m  : \|t\|_2 = 1\}$. It can be shown that $\E g_{\mu} (x) = \nabla f_{\mu}(x)$.
Since we can use only zeroth-order oracle with error we also define the counterpart of the above random gradient-free oracle which can be really computed. We will call it the biased gradient-free oracle:
\begin{equation}
g_{\mu}^{\delta }(x) = \frac{m}{\mu}(f^{\delta}(x+\mu s) -f^{\delta}(x)) s.
\notag
\end{equation}

The following estimates can be proved for the introduced inexact oracle (the full proof is in the Supplementary Materials).
\begin{Lm}
Let $f \in C^{1,1}_{L}$. Then, for any $x,y \in \R^m$,
\begin{align}
& \E \| g_{\mu}^{\delta } (x) \|^2_2 \leq m^2 \mu^2 L^2 + 4 m \|\nabla f(x) \|^2_2 + \frac{8\delta^2 m^2}{\mu^2}\label{expgmd} \\
& - \E \la g_{\mu}^{\delta } (x), x-y \ra \leq- \la \nabla f_{\mu}(x) , x-y \ra + \frac{\delta m}{\mu} \|x-y\|_2.
\label{expgmdxmx}
\end{align}
\label{Lm:gmud}
\end{Lm}

We use gradient-type method with oracle $g_{\mu}^{\delta }(x)$ instead of the real gradient in order to minimize $f_{\mu}(x)$. Since it is uniformly close to $f(x)$ we can obtain a good approximation to the minimum value of $f(x)$.

Algorithm~\ref{alg:GFPGM} below is the variation of the gradient method. Here $\Pi_{X}(x)$ denotes the Euclidean projection of a point $x$ onto a set $X$.
\begin{algorithm}[h!]
        \caption{Gradient-type method}
        \label{alg:GFPGM}
\begin{algorithmic}
   \STATE {\bfseries Input:} The point $x_0$, radius $R$, 
   stepsize $h>0$, number of steps $M$.
         \STATE Define $X=\{ x \in \R^m : \|x-x_0\|_2 \leq 2R \}$.
   \REPEAT
   \STATE Generate $s_k$ and corresponding $g_{\mu}^{\delta }(x_k)$.
   \STATE Calculate $x_{k+1} = \Pi_{X}(x_k- h g_{\mu}^{\delta }(x_k))$.
         \STATE Set $k=k+1$.
   \UNTIL{$k>M$}
         \STATE {\bfseries Output:} The point $x_{k}$.
\end{algorithmic}
\end{algorithm}

Next theorem gives the convergence rate of Algorithm \ref{alg:GFPGM}.
Denote by $\U_k=(s_0, \dots, s_k)$ the history of realizations of the vectors $s_i$, generated on each iteration of the method, $\psi_0 = f(x_0)$, and $\psi_k= \E_{\U_{k-1}}(f(x_{k}))$, $k \geq 1$.

We say that the smooth function is strongly convex with parameter $\tau \geq 0$ if and only if for any $x,y \in \R^m$ it holds that
\begin{equation}
f(x) \geq f(y) + \la \nabla f(y) ,x-y \ra + \frac{\tau}{2} \|x-y\|^2.
\label{eq:fStrConv}
\end{equation}

\begin{Th}
Let $f \in C^{1,1}_{L}$ and the sequence $x_k$ be generated by Algorithm \ref{alg:GFPGM} with $h=\frac{1}{8mL}$.
Then for any $M \geq 0$, we have
\begin{align}
& \frac{1}{M+1}\sum_{i=0}^M \left( \psi_i - f^{*} \right) \leq \frac{8mL R^2}{M+1} + \frac{\mu^2 L (m+8)}{8} + \frac{8 \delta m R}{ \mu }  + \frac{\delta^2 m}{L \mu^2},
\label{eq:rtSmth}
\end{align}
where $f^*$ is the solution of the problem $\min_{x \in \R^m} f(x)$.
If, moreover, $f$ is strongly convex with constant $\tau$, then
\begin{equation}
 \psi_M - f^{*}  \leq \frac12 L \left(\delta_{\mu} + \left(1-\frac{\tau }{16 m L} \right)^M(R^2 - \delta_{\mu})  \right),
\label{eq:rtSmthSC}
\end{equation}
where $\delta_{\mu}=\frac{ \mu^2 L (m+8)}{4 \tau} + \frac{16m \delta R}{\tau \mu } + \frac{2 m \delta^2}{\tau \mu^2 L }  $.
\label{th_1}
\end{Th}


The full of the theorem proof is in the Supplementary Materials. The estimate \reff{eq:rtSmth} also holds for $\hat{\psi}_M \stackrel{\rm def}{=} \E_{\U_{M-1}} f(\hat{x}_M)$, where $\hat{x}_M = \arg \min_x \{ f(x): x \in \{ x_0, \dots, x_M\}\}$.
To make the right hand side of the inequality \reff{eq:rtSmth} less than a desired accuracy $\e$ we need to choose
\begin{align}
& M = \left\lceil \frac{32mLR^2}{\e}\right\rceil, \quad \mu = \sqrt{\frac{2 \e}{L (m+8)}} , \notag \\
& \delta = \min \left\{ \frac{\varepsilon^{\frac32}\sqrt{2}}{32mR\sqrt{L(m+8)}}, \frac{\varepsilon}{\sqrt{2m(m+8)}} \right\} =  \frac{\varepsilon^{\frac32}\sqrt{2}}{32mR\sqrt{L(m+8)}}. \notag
\end{align}

Let's note that we can also estimate the probability of large deviations from the obtained mean rate of convergence. If $f(x)$ is strongly convex, then we have a geometric rate of convergence \reff{eq:rtSmthSC}. Consequently, from the Markov's inequality we obtain that after $O\left(m\frac{L}{\tau}\ln \left(\frac{LR^2}{\e\sigma}\right)\right)$ iterations $\psi_M - f^{*}  \leq \e$ holds with probability greater than $1 - \sigma$. If the function $f(x)$  is not strongly convex, then we can introduce the regularization with parameter $\tau = \e/(2R^2)$ minimizing the function $f(x) + \frac{\tau}{2}\|x\|_2^2$, which is strongly convex. This will give us that after $O\left(m\frac{LR^2}{\e}\ln \left(\frac{LR^2}{\e\sigma}\right)\right)$ iterations $\psi_M - f^{*}  \leq \e$ holds with probability greater than $1 - \sigma$.

\section{Solving the learning problem}
\label{learn}

Our idea for minimizing the function $f(\vp)$ \reff{eq:f_phi_def_2} is the following. We assume that we start from the small vicinity of the optimal value and hence the function $f(\vp)$ is convex in this vicinity (generally speaking, the function \reff{eq:f_phi_def_2} is nonconvex). We choose the desired accuracy $\e$ for approximation of the optimal value of the function $f(\vp)$. This value gives us the number of steps of Algorithm \ref{alg:GFPGM}, the value of the parameter $\mu$, the maximum value of the allowed error of the oracle $\delta$. Knowing the value $\delta$, using Lemma \ref{Lm:Delta_to_delta} we choose the number of steps of the algorithm for an approximate solution of Equation \reff{eq:pi_Phi_P2}, i.e. the number $N$ in \reff{eq:tpi_def}. This idea leads us to Algorithm 2. To the best of our knowledge, this is the first time when the idea of random gradient-free optimization methods is combined with some efficient method for huge-scale optimization using the concept of zero-order oracle with error.

\begin{algorithm}[h!]
        \caption{Method for model learning}
        \label{alg:main_algo}
\begin{algorithmic}
   \STATE {\bfseries Input:} The point $\varphi_0$, $L$ -- Lipschitz constant for the function $f(\varphi)$, radius $R$, 
   accuracy $\varepsilon >0$, numbers $r$, $b$ defined in Lemma \ref{Lm:Delta_to_delta}.
         \STATE Define $X=\{ \vp \in \R^m: \|\vp-\vp_0\|_2 \leq 2R \}$, $M=32 m \frac{LR^2}{\varepsilon}$, $\delta=  \frac{\varepsilon^{\frac32}\sqrt{2}}{32mR\sqrt{L(m+8)}} $, $\mu = \sqrt{\frac{ 2\varepsilon}{L (m+8)}}$.
         \STATE Set $k=0$.
   \REPEAT
   \STATE Generate random vector $s_k$ uniformly distributed over a unit Euclidean sphere $\S$ in $R^m$.
   \STATE Set $N = \frac{1}{\alpha} \ln \frac{ 2 \sqrt{2 r} (2 \sqrt{2 r}  + 2 b)}{\delta}$.
         \STATE For every $q$ from 1 to $|Q|$ calculate  $\tpi^{N}_q(\varphi_k)$, $\tpi^{N}_q(\varphi_k+ \mu s_k)$ defined in \reff{eq:tpi_def}.
         \STATE Calculate  $g_{\mu}^{\delta }(\vp_k) = \frac{m}{\mu}(f^{\delta}(\varphi_k+ \mu s_k)-f^{\delta}(\varphi_k))s_k$, where $f^{\delta}(\varphi)$ is defined in \reff{eq:fdvpDef}.
         \STATE Calculate $\vp_{k+1} = \Pi_{X}\left(\vp_k- \frac{1}{8 m L} g_{\mu}^{\delta }(\vp_k)\right)$.
         \STATE Set $k=k+1$.
         \STATE
   \UNTIL{$k>M$}
         \STATE {\bfseries Output:} The point $\hat{\vp}_M = \arg \min_{\vp} \{ f(\vp): \vp \in \{ \vp_0, \dots, \vp_M\}\}$.
\end{algorithmic}
\end{algorithm}

The most computationally consuming operation on each iteration of the main cycle of this method is the calculation of $2|Q|$ approximate solutions of the equation \reff{eq:pi_Phi_P2}. Hence, each iteration of Algorithm \ref{alg:main_algo} needs approximately $\frac{2|Q|sp}{\alpha} \ln \frac{ 2 \sqrt{2 r} (2 \sqrt{2 r}  + 2 b)}{\delta}$ arithmetic operations, where $s=\max_q{s_q}$, $p=\max_q{p_q}$.
So, we obtain the following theorem, which gives the result for local convergence of Algoritghm \ref{alg:main_algo}.
\begin{Th}
Assume that the point $\varphi_0$ lies in the vicinity of the local minimum point $\varphi^*$ of the function $f(\varphi)$ and the function $f(\varphi)$ is convex in this vicinity. Then the mean total number of arithmetic operations for the accuracy $\e$ (i.e. for inequality $\E_{\U_{M-1}} f(\hat{\vp}_M) - f(\vp^*) \leq \e$ to hold) is given by
\begin{align}
&64 mps|Q| \frac{LR^2}{\alpha \varepsilon}  \ln \left(4 (2 r  + b\sqrt{2 r})  \frac{32mR\sqrt{L(m+8)} }{\varepsilon^{\frac32}\sqrt{2}} \right). \notag
\end {align}
\label{th_2}
\end{Th}


Let us make some remarks. Note that each iteration of the main cycle of the algorithm above can be fully paralleled using $|Q|$ processors. Also it is important that the use of geometrically convergent method as the inner algorithm leads to the overall complexity bound which is the product of complexity bounds of the inner and outer algorithms.


The direct calculation of the parameter $L$ has many obstacles and leads to the overestimation. Another way is to use the restart method. Since we know the exact required number of iterations for the fixed accuracy, confidence level and $L$, we can use the following procedure. We start with some initial value of $L$. Calculate the approximation by Algorithm \ref{alg:main_algo}. Then set $L := 2L$ and repeat, i.e. calculate the approximation by the Algorithm \ref{alg:main_algo}, working with new $L$, etc. The stopping criterion here is stabilization (with the same accuracy as before) of this sequence of function values. The total number of such restarts will be of the order $\log_2(2L)$. The same can be done with the unknown parameter $R$.

Here we have omitted the full description of the generalization of the fast-gradient-type scheme~[14,~16] for the case of inexact oracle and application of the obtained method for the minimization of the function $f(\vp)$. The fast-gradient-type scheme is faster but requires the oracle to be more precise. The resulting mean value of the number of arithmetic operations to achieve the accuracy $\e$ for this method is
\begin{align}
& O \left( mps|Q| \sqrt{\frac{LR^2}{\alpha^2 \varepsilon}}  \ln \left((r  + b\sqrt{r})  \frac{mRL}{\e}   \right) \right). \notag
\end {align}

\begin{algorithm}[tb]
        \caption{Fast method for model learning}
        \label{alg:main_algo_fast}
\begin{algorithmic}
   \STATE {\bfseries Input:} The point $\varphi_0$, $L$ -- Lipschitz constant for the function $f(\varphi)$, $\tau$ -- the strong convexity parameter of the function $f(\varphi)$ (note that $\tau=0$ if the function is convex), number $R$ such that $\|\varphi_0-\varphi^*\|_2 \leq R$,  accuracy $\varepsilon >0$, numbers $r$, $b$ defined in Lemma \ref{Lm:Delta_to_delta}.
         \STATE Define $N=16 m \sqrt{\frac{3LR^2}{\varepsilon}}$, $\mu = \sqrt{\frac{ 64 \varepsilon}{3L (5N+64)}}$, $\delta= \sqrt{ \frac{4 \varepsilon \mu^2 L }{3 N} }$, $\gamma_0=L$, $v_0 = \vp_0$, $\theta = \frac{1}{64m^2L}$, $h=\frac{1}{8mL}$.
         \STATE Set $k=0$.
   \REPEAT
   \STATE Compute $\alpha_k > 0$ satisfying $\frac{\alpha_k^2}{\theta} = (1-\alpha_k) \gamma_k + \alpha_k \tau \equiv \gamma_{k+1}$.
   \STATE Set $\lambda_k = \frac{\alpha_k}{\gamma_{k+1}}\tau$, $\beta_k = \frac{\alpha_k \gamma_k }{\gamma_k + \alpha_k \tau }$, and $ y_k = (1-\beta_k)\vp_k + \beta_k v_k$.
         \STATE Generate random vector $s_k$ uniformly distributed over a unit Euclidean sphere $\S$ in $R^m$
         \STATE Set $\hat{N} = \frac{1}{\alpha} \ln \frac{ 2 \sqrt{2 r} (2 \sqrt{2 r}  + 2 b)}{\delta}$.
         \STATE For every $q$ calculate  $\tpi^{\hat{N}}_q(\varphi_k)$, $\tpi^{\hat{N}}_q(\varphi_k+ \mu s_k)$ defined in \reff{eq:tpi_def}.
         \STATE Calculate  $g_{\mu}^{\delta }(\vp_k) = \frac{m}{\mu}(f_{\delta}(\varphi_k+ \mu s_k)-f_{\delta}(\varphi_k))s_k$, where $f_{\delta}(\varphi)$ is defined in \reff{eq:fdvpDef}.
         \STATE Calculate $\vp_{k+1} = y_k - h g_{\mu}^{\delta }(y_k)$, $ v_{k+1}=(1- \lambda_k) v_k + \lambda_k y_k - \frac{\theta}{\alpha_k}g_{\mu}^{\delta }(y_k)$.
         \STATE Set $k=k+1$.
         \STATE
   \UNTIL{$k>N$}
         \STATE {\bfseries Output:} The point $\vp_N$.
\end{algorithmic}
\end{algorithm}


Also we want to point that the algorithm for solving equation \reff{eq:pi_Phi_P2} was chosen consciously from a set of modern methods for computing PageRank. We used review~[5] of such methods. Since for our problem we need to estimate the error which is introduced to the function $f(\vp)$ value by approximate solution of the ranking problem \reff{eq:pi_Phi_P2}, we considered only three methods: Markov Chain Monte Carlo (MCMC), Spillman's and Nemirovski-Nesterov's (NN). These three methods allow to make the difference $\|\pi_q(\vp)-\tilde{\pi}_q\|$, where $\tilde{\pi}_q$ is the approximation, small. This is crucial to prove results like Lemma \ref{Lm:Delta_to_delta}. Spillman's alogoritm converges in infinity norm which is usually $\sqrt{p}$ times larger than 2-norm. MCMC converges in 2-norm and NN converges in 1-norm. Finally, the full complexity analysis of the two-level algorithm showed that for the dimensions $m,p$ and accuracy $\e$ considered in our work the combination of gradient-free method with NN method is better than the combination with MCMC in terms of upper bound for arithmetic operations needed to achieve given accuracy.

\section{Experimental results}
\label{Experimental results}

We compare the performances of different learning techniques, our gradient-free method, an untuned gradient-free method and classical PageRank. In the next section, we describe the graph, which we exploit in our experiments (the user browsing graph). In Section~\ref{data} and Section~\ref{results}, we describe the dataset and the results of the experiments respectively.

\subsection{User browsing graph}
\label{graph}

In this section, we define the web user browsing graph (which was first considered in~[12]). We choose the user browsing graph instead of a link graph with the purpose to make the model query-dependent.

Let $q$ be any query from the set $Q$. A user session $S_q$ (see~[12]), which is started from $q$, is a sequence of pages $(i_1,i_2,...,i_k)$ such that, for each $j\in\{1,2,...,k-1\}$, the element $i_j$ is a web page and there is a record $i_j\rightarrow i_{j+1}$ which is made by toolbar. The session finishes if the user types a new query or if more than 30 minutes left from the time of the last user's activity. We call pages $i_j,i_{j+1}$, $j\in\{1,\ldots,k-1\}$, {\it the neighboring elements} of the session $S_q$.

We define the user browsing graph $\Gamma=(V,E)$ as follows. The set of vertices $V$ consists of all the distinct elements from all the sessions which are started from any query $q\in Q$. The set of directed edges $E$ represents all the ordered pairs of neighboring elements $(\tilde i,i)$ from the sessions. For any $q\in Q$, we set $F_q(\varphi_1,i)=0$ for all $\varphi_1\in\mathbb{R}^{m_1}$ if there is no session which is started from $q$ and contains $i$ as its first element. Moreover, we set $G_q(\varphi_2,\tilde i\rightarrow i)=0$ for all $\varphi_2\in\mathbb{R}^{m_2}$ if there is no session which is started from $q$ and contains the pair of neighboring elements $\tilde i,i$.

As in~[21], we suppose that for any $q\in Q$ , any $i\in V^1_q$ and any $\tilde i\rightarrow i\in E_q$, a vector of node's features $\mathbf{V}^q_i\in\mathbb{R}^{m_1}$ and a vector of edge's features $\mathbf{E}^q_{\tilde i i}\in\mathbb{R}^{m_2}$ are given. We set $F_q(\varphi_1,i)=\langle\varphi_1,\mathbf{V}^q_i\rangle$, $G_q(\varphi_1,\tilde i\rightarrow i)=\langle\varphi_2,\mathbf{E}^q_{\tilde i i}\rangle$.

\subsection{Data}
\label{data}

All experiments are performed with pages and links crawled by a popular commercial search engine. We utilize all the records from the toolbar that were made from 27 October 2014 to 18 January 2015. We randomly choose the set of queries $Q$ the user sessions start from, which contains $\approx1$K queries. There are $\approx0.6$M vertices and $\approx0.8$M edges in graphs $\Gamma_q$, $q\in Q$, in total. For each query a set of pages was judged by professional assessors hired by the search engine. Our data contains $\approx3.8$K judged query--document pairs. The relevance score is selected from among 5 editorial labels. We divide our data into two parts. On the first part ($80\%$ of the set of queries $Q$) we train the parameters and on the second part we test the algorithms.
To define weights of nodes and edges we consider a set of 26 query--document features.
For any $q\in Q$ and $i\in V^1_q$, the vector $\mathbf{V}^q_i$ contains values of all these feautures for query--document pair $(q,i)$. We set $m_2=2m_1=52$ and $\mathbf{E}^q_{\tilde i,i}=([\mathbf{V}^q_{\tilde i}]_1,\ldots,[\mathbf{V}^q_{\tilde i}]_{m_1},[\mathbf{V}^q_i]_1,\ldots,[\mathbf{V}^q_i]_{m_1})$.

\subsection{Ranking quality}
\label{results}

We find the optimal values of the parameters for all the methods by minimizing the objective $f$ defined by Equation~\ref{eq:f_phi_def_1} by the common untuned gradient-free method GF1 (Algorithm~\ref{alg:GFPGM}) and our precise gradient-free method GF2 (Algorithm~\ref{alg:main_algo}). Besides, we use PageRank (PR) as the common baseline for the algorithms (used as the only baseline for SSP in~[6] and one of the baselines for SNP in~[21]).

The sets of parameters which are exploited by the optimization methods (and not tuned by them) are the following: the Lipschitz constant $L=1.6\cdot 10^{-4}$, the accuracy $\varepsilon=6.9\cdot 10^{-3}$ (in GF2), the radius $R=1$ (in both GF1 and GF2), the parameter $N=117$~(\ref{eq:tpi_def}), which defines the approximation $\tilde\pi^N_q$ of the stationary distribution $\tilde\pi_q$, of algorithms GF1 and PR is chosen in such a way that the accuracy $\Delta$~(\ref{Delta}) equals $10^{-8}$. Moreover, $M=10$ (the number of iterations of the optimization method) and $h=10$ (the stepsize) in the algorithms GF1 (the number of iterations is less than the value of this parameter in GF2).

In Table~\ref{Table:NDCG}, we present the ranking performances in terms of our loss function $f$.

\begin{table}[!ht]
\centering
\begin{tabular}{|c|c|}

  \hline

Method         &  f~(Equation~\ref{eq:f_phi_def_1}) \\

  \hline

  \hline

  GF2                      &    0.00107   \\
  
  \hline

  GF1                     &     0.001305   \\
  
  \hline

  

  PR                     &      0.0118   \\
  
  \hline  
  
\end{tabular}
\vspace{-0.3cm}

\caption{\small Performances of GF2, GF and PR methods.} \label{Table:NDCG}
\end{table}

\vspace{-0.3cm}

Moreover, the NDCG@3 (@5) gains of both GF1 and GF2 in comparison with PR exceeds $20\%$ for both metrics. We obtain the $p$-values of the paired $t$-tests for all the above differences in ranking qualities on the test set of queries. These values are less than 0.005. Thus, we conclude that the obtained values of the parameters by our optimization method are closer to optimal than in the case of GF1.

\section{Conclusion}
\label{Conclusion}


We consider a problem of learning parameters of supervised PageRank models, which are based on calculating the stationary distributions of the Markov random walks with transition probabilities depending on the parameters. Due to huge hidden dimension of the optimization problem and the impossibility of exact calculating derivatives of the stationary distributions w.r.t. its parameters, we propose a two-level method, based on random gradient-free method with inexact oracle to solve it instead of the previous gradient-based approach. We find the best settings of the gradient-free optimization method in terms of the number of arithmetic operations needed to achieve given accuracy of the objective. In particular, for the proposed method, we provide an estimate for the total number of arithmetic operations to obtain the given accuracy in terms of local convergence. We apply our algorithm to the web page ranking problem by considering a dicrete-time Markov random walk on the user browsing graph. Our experiments show that our two-level method outperforms both classical PageRank algorithm and the gradient-free algorithm with other settings (which are, theoretically, not optimal). In the future, some globalization techniques can be considered (e.g., multi-start), because the objective function is nonconvex.


\subsubsection*{Acknowledgment}
The work was partially supported by Russian Foundation for Basic Research grants 14-01-00722-a, 15-31-20571-mol\_a\_ved.

\subsubsection*{References}

\small{
[1] A. Agarwal, O. Dekel, L. Xiao, {\it Optimal algorithms for online convex optimization with multi-point bandit feedback}, COLT'2010.

[2] L.~Backstrom, J.~Leskovec, {\it Supervised random walks: predicting and recommending links in
social networks}, WSDM'11.

[3] Na~Dai. Brian~D.~Davison, {\it Freshness Matters: In Flowers, Food, and Web Authority}, SIGIR'10.


[4] N.~Eiron, K.~S.~McCurley, J.~A.~Tomlin, {\it Ranking the web frontier}, WWW'04.

[5] A. Gasnikov, D. Dmitriev, {\it  Efficient randomized algorithms for PageRank problem},  Comp. Math. \& Math. Phys, 2015, V.~55, No.~3, P.~1--18. 

[6] B.~Gao, T.-Y.~Liu, W.~W.~Huazhong, T.~Wang, H.~Li, {\it Semi-supervised ranking on very large graphs with rich metadata}, KDD'11.

[7] S. Ghadimi, G. Lan, {\it Stochastic first- and zeroth-order methods for nonconvex stochastic programming}, SIAM Journal on Optimization, 2014, 23(4), 2341–-2368.

[8] T.~H.~Haveliwala, {\it Efficient computation of PageRank}, Stanford University Technical Report, 1999.

[9] T.~H.~Haveliwala, {\it  Topic-Sensitive PageRank}, WWW'02.

[10] G.~Jeh, J.~Widom, {\it Scaling Personalized Web Search}, WWW'03.

[11] J.~M.~Kleinberg, {\it Authoritative sources in a hyperlinked environment}, SODA'98.

[12] Y.~Liu, B.~Gao, T.-Y.~Liu, Y.~Zhang, Z.~Ma, S.~He, H.~Li, {\it BrowseRank: Letting Web Users Vote for Page Importance}, SIGIR'08.

[13] J.~Matyas, {\it Random optimization}, Automation and Remote Control, 1965, V.~26, P.~246-253.

[14] Yu.~Nesterov, {\it Introductory Lectures on Convex Optimization}, Kluwer, Boston, 2004.

[15] Yu.~Nesterov, {\it Efficiency  of  coordinate  descent  methods  on  huge-scale optimization problems}, SIAM Journal on Optimization, 2012, V. 22, No.2, p. 341--362.

[16] Yu.~Nesterov, {\it Random gradient-free minimization of convex functions}, CORE Discussion Paper, 2011/1 \url{http://www.uclouvain.be/cps/ucl/doc/core/documents/coredp2011_1web.pdf}.

[17] Yu.~Nesterov, A.~Nemirovski, Finding the stationary states of Markov chains by iterative methods {\it CORE Discussion Paper 2012/58}, Applied Mathematics and Computation, (2014), \url{http://dial.academielouvain.be/handle/boreal:122163}.

[18] L.~Page, S.~Brin, R.~Motwani, and T.~Winograd, {\it The PageRank citation ranking: Bringing order to the web}, \url{http://dbpubs.stanford.edu/pub/1999-66}, 1999.

[19] M.~Richardson, P.~Domingos, {\it The intelligent surfer: Probabilistic combination of link and
content information in PageRank}, NIPS'02.


[20] M.~Zhukovskii, G.~Gusev, P.~Serdyukov, {\it URL Redirection Accounting for Improving Link-Based Ranking Methods}, ECIR'13.

[21] M.~Zhukovskiy, G.~Gusev, P.~Serdyukov, {\it Supervised Nested PageRank}, CIKM'14.

[22] M.~Zhukovskii, A.~Khropov, G.~Gusev, P.~Serdyukov, {\it Fresh BrowseRank}, SIGIR'13.


}

\newpage

\section{Appendix}
\subsection{Proof of Lemma \ref{Lm:gmud}}

We will need the following lemma.

\begin{Lm}
Let $s$ be random vector uniformly distributed over the unit sphere $\S \in \R^m$. Then
\begin{equation}
\E_s(\la \nabla f(x), s \ra)^2 = \frac{1}{m}\|\nabla f(x)\|^2_*.
\label{expnfss}
\end{equation}
\end{Lm}

\myproof
We have $\E_s(\la \nabla f(x), s \ra)^2 = \frac{1}{S_m(1)} \int_{S^m} (\la \nabla f(x), s \ra)^2 d \sigma(s)$, where $S_m(r)$ is the volume of the unit sphere which is the border of the ball in $\R^m$ with radius $r$. Note that $S_m(r)=S_m(1) r^{m-1}$.
Let $\varphi$ be the angle between $\nabla f(x)$ and $s$. Then
\begin{align}
&\frac{1}{S_m(1)} \int_{S^m} (\la \nabla f(x), s \ra)^2 d \sigma(s) = \frac{1}{S_m(1)}  \int_{0}^{\pi} \|\nabla f(x)\|^2_* \cos^2 \varphi S_{m-1}( \sin \varphi) d \varphi = \notag \\
&= \frac{S_{m-1}(1)}{S_m(1)}  \|\nabla f(x)\|^2_* \int_{0}^{\pi}\cos^2 \varphi \sin^{m-2} \varphi d \varphi
\notag
\end{align}
First changing the variable using equation $x=\cos \varphi$, and then $t=x^2$, we obtain
$$
\int_{0}^{\pi}\cos^2 \varphi \sin^{m-2} \varphi d \varphi = \int_{-1}^{1}x^2 (1-x^2)^{(m-3)/2}d x =  \int_{0}^{1}t^{1/2} (1-t)^{(m-3)/2}d t =  B\left(\frac32,\frac{m-1}{2} \right)=  \frac{\sqrt{\pi} \Gamma\left(\frac{m-1}{2}\right)}{2 \Gamma\left(\frac{m+2}{2}\right)},
$$
where $\Gamma(\cdot)$ is the Gamma-function.
Also we have
\begin{equation}
\frac{S_{m-1}(1)}{S_m(1)} = \frac{m-1}{m \sqrt{\pi}} \frac{\Gamma\left(\frac{m+2}{2}\right)}{\Gamma\left(\frac{m+1}{2}\right)}.
\label{eq:snm1tosn}
\end{equation}

Finally using the relation $\Gamma(m+1) = m \Gamma(m)$, we obtain
$$
\E(\la \nabla f(x), s \ra)^2 = \|\nabla f(x)\|^2_* \left(1-\frac{1}{m}\right) \frac{ \Gamma\left(\frac{m-1}{2}\right)}{ 2\Gamma\left(\frac{m+1}{2}\right)}=\|\nabla f(x)\|^2_* \left(1-\frac{1}{m}\right) \frac{ \Gamma\left(\frac{m-1}{2}\right)}{2 \frac{m-1}{2}\Gamma\left(\frac{m-1}{2}\right)} = \frac{1}{m}\|\nabla f(x)\|^2_*
$$
\qed


Using \reff{eq:fLipSm} we obtain
\begin{align}
& (f_{\delta}(x+\mu s) - f_{\delta}(x))^2 = \notag \\
& (f(x+\mu s) - f(x) - \mu \la \nabla f(x), s \ra + \mu \la \nabla f(x), s \ra + \tilde{\delta}(x+\mu s) - \tilde{\delta}(x))^2 \leq \notag \\
& 2 (f(x+\mu s) - f(x) - \mu \la \nabla f(x), s \ra + \mu \la \nabla f(x), s \ra)^2 + 2(\tilde{\delta}(x+\mu s) - \tilde{\delta}(x))^2 \leq \notag \\
& 4 \left(\frac{\mu^2}{2}L_1 \|s\|^2\right)^2 + 4 \mu^2 (\la \nabla f(x), s \ra)^2 + 8 \delta^2 = \mu^4 L_1^2 \|s\|^4 + 4 \mu^2 (\la \nabla f(x), s \ra)^2 + 8 \delta^2 \notag
\end{align}
Using \reff{expnfss}, we get
\begin{align}
& \E_s \| g_{\mu}^{\delta } (x) \|^2_* \leq \frac{m^2}{\mu^2 V_s} \int_S \left(\mu^4 L_1^2 \|s\|^4 + 4 \mu^2 (\la \nabla f(x), s \ra)^2 + 8 \delta^2 \right) \|s\|^2_* d\sigma (s) = m^2 \mu^2 L_1^2 + 4 m \|\nabla f(x) \|^2_* + \frac{8\delta^2 m^2}{\mu^2}.\notag
\end{align}

Using the equality $\E_sg_{\mu (x)} = \nabla f_{\mu}(x)$, we have
\begin{align}
& - \E_s \la g_{\mu}^{\delta } (x), x-x^* \ra = - \frac{m}{\mu V_s} \int_S (f_{\delta}(x+\mu s) - f_{\delta}(x)) \la s, x-y \ra  d\sigma (s) = \notag \\
& = - \frac{m}{\mu V_s} \int_S (f(x+\mu s) - f(x)) \la s, x-y \ra  d\sigma (s) - \notag \\
& - \frac{m}{\mu V_s} \int_S (\tilde{\delta}(x+\mu s) - \tilde{\delta}(x)) \la s, x-y \ra  d\sigma (s) \leq - \la \nabla f_{\mu}(x) , x-y \ra + \frac{\delta m}{\mu} \|x-y\|.
\notag
\end{align}

\qed

\subsection{Proof of Theorem \ref{th_1}}

We extend the proof in~[16] for the case of randomization on a sphere (instead of randomization based on normal distribution) and for the case when one can calculate the function value only with some error of unknown nature.

Consider the point $x_k$, $k\geq 0$ generated by the method on the $k$-th iteration. Denote $r_k=\|x_k-x^*\|_2$. Note that $r_k \leq 4R$. We have:
\begin{align}
&r_{k+1}^2 = \|x_{k+1}-x^*\|_2^2 \leq  \|x_{k}-x^*-h g_{\mu}^{\delta } (x_k)\|_2^2 = \notag \\
& = \|x_{k}-x^*\|_2^2 - 2h \la g_{\mu}^{\delta } (x_k), x_k-x^* \ra + h^2 \|g_{\mu}^{\delta } (x_k) \|_2^2.
\notag
\end{align}

Taking the expectation with respect to $s_k$ we get
\begin{align}
& \E_{s_k} r_{k+1}^2 \stackrel{\reff{expgmd},\reff{expgmdxmx}}{\leq}  r_k^2 - 2h\la \nabla f_{\mu}(x_k) , x_k-x^* \ra + \frac{2 \delta m h }{\mu} r_k + \notag \\
& + h^2 \left( m^2 \mu^2 L^2 + 4 m \|\nabla f(x_k) \|_2^2 + \frac{8\delta^2 m^2}{\mu^2} \right) \leq \notag \\
& \leq  r_k^2 -2h (f(x_k )- f_{\mu}(x^*)) + \frac{8 \delta m h R}{\mu} + \notag \\
& + h^2 \left( m^2 \mu^2 L^2 + 8 m L (f(x_k) - f^*) + \frac{8\delta^2 m^2}{\mu^2} \right) \leq \notag \\
& \leq r_k^2 -2h (1-4h m L) (f(x_k )- f^*) + \frac{8 \delta m h R}{\mu}  +  \notag \\
&+ m^2 h^2 \mu^2 L^2 + h L \mu^2 + \frac{8\delta^2 m^2 h^2}{\mu^2} \leq \notag \\
& \leq r_k^2 + \frac{R \delta  }{ \mu L} - \frac{f(x_k )- f^*}{8mL} + \frac{\mu^2(m+8)}{64 m} + \frac{\delta^2 }{8\mu^2L^2}.
\label{ThalgSmthPr1}
\end{align}
Taking expectation with respect to $\U_{k-1}$ and defining $\rho_{k+1} \stackrel{\rm def}{=}\E_{\U_{k}} r_{k+1}^2$ we obtain
\begin{equation}
\rho_{k+1} \leq \rho_k- \frac{\psi_k- f^*}{8mL}+ \frac{\mu^2(m+8)}{64 m} + \frac{R \delta  }{ \mu L} +  \frac{\delta^2 }{8\mu^2L^2}.
\notag
\end{equation}
Summing up these inequalities and dividing by $N+1$ we obtain \reff{eq:rtSmth}.

Now assume that the function $f(x)$ is strongly convex. From \reff{ThalgSmthPr1} we get
\begin{equation}
\E_{s_k} r_{k+1}^2 \stackrel{\reff{eq:fStrConv}}{\leq} \left(1-\frac{\tau}{16mL} \right) r_k^2 + \frac{R \delta  }{ \mu L} + \frac{\mu^2(m+8)}{64 m} + \frac{\delta^2 }{8\mu^2L^2}
\notag
\end{equation}
Taking expectation with respect to $\U_{k-1}$ we obtain
\begin{equation}
\rho_{k+1} \leq \left(1-\frac{\tau}{16mL} \right) \rho_k + \frac{R \delta  }{ \mu L} + \frac{\mu^2(m+8)}{64 m} + \frac{\delta^2 }{8\mu^2L^2}
\notag
\end{equation}
and
\begin{align}
& \rho_{k+1} - \delta_{\mu} \leq \left(1-\frac{\tau}{16mL} \right) ( \rho_k - \delta_{\mu}) \leq \notag \\
& \leq \left(1-\frac{\tau}{16mL} \right)^{k+1} ( \rho_0 - \delta_{\mu}). \notag
\end{align}

Using the fact that $\rho_0 = R^2$ and $\psi_k - f^* \leq \frac12 L \rho_k$ we obtain \reff{eq:rtSmthSC}.
\qed

\end{document}